\documentclass[12pt,reqno]{amsart}

\usepackage{amsfonts,amssymb,latexsym,amsmath}

\usepackage{xcolor}
\usepackage{amssymb}
\usepackage{tikz}
\usepackage{tabularx}

\newtheorem{theorem}{Theorem}[section]
\newtheorem{lemma}[theorem]{Lemma}

\newtheorem{proposition}[theorem]{Proposition}

\theoremstyle{definition}

\numberwithin{equation}{section}

\DeclareMathOperator{\tr}{trace}

\newcommand{\T}{{\mathbb T}}
\renewcommand{\kappa}{\varkappa}
\newcommand{\be}{\begin{equation}}
\newcommand{\ee}{\end{equation}}

\newcommand{\bq}{\begin{eqnarray}}
\newcommand{\eq}{\end{eqnarray}}
\newcommand{\nn}{\nonumber}
\newcommand{\ba}{\begin{array}}
\newcommand{\ea}{\end{array}}

\newcommand{\iv}{^{-1}}
\newcommand{\iy}{\infty}
\newcommand{\al}{\alpha}

\newcommand{\cL}{\mathcal{L}}

\newcommand{\cS}{\mathcal{S}}

\title{Factoring determinants and applications to number theory}
\author{Estelle Basor}
\address{American Institute of Mathematics, Pasadena, CA, USA}
\email{ebasor@aimath.org}
\author{Brian Conrey      }
\address{American Institute of Mathematics, Pasadena, CA, USA}
\email{conrey@aimath.org}

\begin{document}

\maketitle

 \begin{abstract}  Products of shifted characteristic polynomials, and ratios of such products, averaged over the classical compact groups are of great interest to number theorists as they model similar averages of L-functions in families with the same symmetry type as the compact group.
  We use Toeplitz and Toeplitz plus Hankel operators and the identities
  of Borodin - Okounkov - Case - Geronimo, and Basor - Ehrhardt to prove that, in certain cases, these unitary averages factor as polynomials into averages over the symplectic group and the orthogonal group. Building on these identities we present  new proofs of the exact formulas for these averages where the ``swap'' terms that are 
  characteristic of the number theoretic averages occur from the  Fredholm expansions of the determinants of the appropriate Hankel operator.  This is the fourth different proof of the formula for the averages of ratios of products of shifted characteristic polynomials; the other proofs are based on supersymmetry; symmetric function theory, and orthogonal polynomial methods from Random Matrix Theory. 
\end{abstract}

\section{introduction}

The statistical study of L-functions is the backbone of modern analytic number theory. The quantities of primary interest are the zeros and their distributions and the values 
both large and small. One approaches these statistics by calculating averages (moments)  of
 products of shifted L-functions and averages of ratios of such. The averages are over a family of L-functions 
 with conductor up to some bound.

Our understanding of moments and ratios of L-functions in families has increased dramatically in the last 25 years
stimulated initially by the work of Keating and Snaith [KeSn] showing that Random Matrix Theory provides
 a good model for these moments and by the work of Katz and Sarnak [KaSa] that showed how to associate a symmetry type to each family of L-functions.
 The symmetry type is one of Unitary, Symplectic or Orthogonal and governs which 
 classical group with Haar measure one should use in trying to model a given family of 
 L-functions.
 
Two conjectures about these averages have come into prominence: the {\it} recipe tells how to calculate moments in a family, and the ratios conjecture
tells how to average ratios of products of L-functions. Each of these has a corresponding theorem in Random Matrix Theory (RMT), (see [CFKRS, I and II], [BG], [CFZ], [CFS].)
RMT has become an indispensable tool for predicting the answers for questions about L-functions.
For example, if one wants to know the distribution of the real parts of the zeros of $\zeta'(s)$ one would might start by asking  the same question about the radial distribution of zeros of $\Lambda_U'(z)$
where $U\in U(N)$, the group of $N\times N$ unitary matrices, and $\Lambda_U(x)=\det(I-xU^* )$ is the characteristic polynomial of $U$. If one could evaluate
$$\int_{U(N)} |\Lambda_U'(x)|^{2k}dU $$ 
for small positive  real $k$ and $x$ near 1 one could answer this (unsolved) problem and the answer, we believe, would be the same distribution 
(with appropriate scaling) as the question for zeros of $\zeta'(s)$.

In general, moments of L-functions give information about the value distribution (in particular large values) of L-functions, while averages of  ratios
of products of L-functions give us information about the distribution of small values, in particular zeros, of L-functions (see [CFKRS], [CS07], and  [CS08] for some examples).

Examples of families of L-functions include:
\begin{itemize}
\item $\{\zeta(1/2+i t):t\in \mathbb R\}$ is a unitary family.
In fact any L-function regarded in $t$ aspect is a family parametrized by the (continuous variable) $t$. All other families are discrete.
\item $\{L(1/2,\chi_d): \chi_d \mbox{ is a real primitive Dirichlet character}\}$ is a symplectic family. These characters are extensions of the Legendre symbol $\chi_p(n)=\big(\frac n p\big)$ which is 1 if $n$ is a square modulo $p$ and $-1$ if $n$ is not a square.
\item $\{L_f(1/2):f\mbox{ is a weight 2, level $q$ newform}\}$ is an orthogonal family. Examples include the L-functions associated with elliptic curves.  This is a family of degree 2 L-functions.
\end{itemize}

The symbiotic relationship between Analytic Number Theory (ANT) and RMT reaps many rewards. Conjectures we make in ANT have counterparts that are expected to be proveable in RMT. These are often of a different nature than random matrix theorists might typically think about.  And theorems in RMT lead us to wonder whether the analogue in ANT is true and so to conjectures in ANT that guide the way as we try to understand  L-functions.

In this paper we prove a theorem in RMT that would be very surprising in ANT. 
It has to do with moments. For example we might consider
$$ \int_0^\infty \psi(\frac t T) \zeta(1/2+it + \alpha)\zeta(1/2-it +\beta) ~dt$$ where $\psi $ is (a) compactly supported (test-function) on $[1,2].$
The RMT analogue of this simplest shifted moment problem would be 
$$\int_{U(N)} \det(I-a~U) \det(I-b~U^*)dU$$
where we integrate with respect to Haar measure $dU$ and where $a=e^{-\alpha}$. In ANT we can prove that
\begin{eqnarray*}&&
 \int_0^\infty \psi(\frac t T) \zeta(1/2+it + \alpha)\zeta(1/2-it +\beta) ~dt\\&&\quad = \int_0^\infty \psi(\frac t T) \bigg(\zeta(1+\alpha+\beta)+\big(\frac t{2\pi}\big)^{-\alpha-\beta}\zeta(1-\alpha-\beta)\bigg)~dt +O(T^{\frac 1 3 })
\end{eqnarray*}
for $\alpha, \beta \ll (\log T)^{-1}$
and in RMT we can prove that
$$\int_{U(N)} \det(I-e^{-\alpha} U)\det(I-e^{-\beta}U^*)~dU = z(\alpha+\beta)+e^{-N(\alpha+\beta)} z(-\alpha-\beta)$$
for any $\alpha$ and $\beta$ where $$z(x)=\frac{1}{1-e^{-x}}.$$

The surprising theorem that we prove is that in certain cases the moments of a product of shifted characteristic polynomials over the unitary group 
factors into a moment over the symplectic group times a moment over the orthogonal group. We do not expect  this to be literally true 
for moments of families of L-functions but it is intriguing to try to speculate what an analogue might look like.
 
There are examples of similar  phenomena that have been previously observed in RMT. For example,  see  section 7.6 of [KaSa] and the many formulas in [BR].

This factorization also works for unitary averages of ratios of products of shifted characteristic polynomials.

The proof depends on some remarkable identities. The first is known as the Borodin-Okounkov-Case-Geronimo (BOCG) identity which holds for for unitary averages. The second is a non-trivial  extension to symplectic and orthogonal averages discovered 
and proven by Basor and Ehrhardt. 
 
We give some examples of factorizations:
We have 
\begin{eqnarray*}&& \int_{U(3)} \Lambda_U(a) \Lambda_U(b) \Lambda_{U^*}(a) \Lambda_{U^*}(b)~dU\\ &&=(1 + a^2 + a b + b^2 + a^2 b^2)\\&&\qquad \times (1 + a^4 + a b + a^3 b + 2 a^2 b^2 + 
   a b^3 + a^3 b^3 + b^4 + a^4 b^4)\end{eqnarray*}
 while
   \begin{eqnarray*}
   \int_{USp(2)}\Lambda_U(a)\Lambda_U(b)dU=1 + a^2 + a b + b^2 + a^2 b^2
   \end{eqnarray*}
   and
 \begin{eqnarray*}\int_{O^{+}(4)}\Lambda_U(a)\Lambda_U(b)dU&=&1 + a^4 + a b + a^3 b + 2 a^2 b^2 \\&&\qquad + a b^3 + a^3 b^3 + b^4 + a^4 b^4.\end{eqnarray*}

As an example of a  ratios calculations we mention
 \begin{eqnarray*} &&\int_{U(N)}\frac{  \Lambda_U(a)  \Lambda_{U^*}(a)}
{ \Lambda_U(c)  \Lambda_{U^*}(c)}dU \\&&\qquad =\frac{(1 + a^{N+1} - a c - a^N c) (1 - a^{N+1} - a c + a^N c)}{(1 - a) (1 + a) (1 - c) (1 + c)}
\end{eqnarray*}
while for even $N$
 \begin{eqnarray*}
   \int_{USp(N)}\frac{\Lambda_U(a)}{\Lambda_U(c)}~dU= \frac{(1 - a^{N+2} - a c + a^{N+1} c)  }{(1 - a) (1 + a)  }
   \end{eqnarray*}
   and
\begin{eqnarray*}
 \int_{O^(+N)}\frac{\Lambda_U(a)}{\Lambda_U(c)}~dU
 &=& \frac{(1 + a^{N} - a c - a^{N-1} c)  }{(1 - c) (1 + c)  }.
 \end{eqnarray*}
Thus, for  $N$ odd, we have 
\begin{eqnarray*}
\int_{U(N)}\frac{  \Lambda_U(a)  \Lambda_{U^*}(a)}
{ \Lambda_U(c)  \Lambda_{U^*}(c)}dU= \left( \int_{USp(N-1)}\frac{\Lambda_U(a)}{\Lambda_U(c)}dU\right)\left(\int_{O^+(N+1)}\frac{\Lambda_U(a)}{\Lambda_U(c)}~dU\right).
\end{eqnarray*}

 In the process of recording  the proof of this factorization, we discovered a new proof for the ratios theorem. This is the fourth known proof. The first was [CFZ] motivated by physics and uses
  supersymmetry; the second was [BG] using symmetric function theory and representation theory. The third was [CFS] using orthogonal polynomial methods from RMT.
 The new one uses operator theory and has the nice feature that the idea is very simple to explain. The unitary moment can be expressed via BOCG as a manifestation of the Szego-Widom theorem for Toeplitz matrices. One piece of this identity involves a finite rank  operator $K$ whose Fredholm expansion  gives the desired swap terms.
 To handle the formulas for symplectic and orthogonal one proceeds in the same way using the Basor-Ehrhardt identity which relates to operators that are ``Toeplitz plus Hankel.''
The first proof works for all matrix size $N\ge 1$, whereas the second and third proofs have a condition that $N$ must be bigger than the number of factors in the average. 
Our new proof works for all $N\ge 1$.

Work of the second author was supported in part by a Focused Research Grant DMS 1854398 from the NSF.

\section{Paper organization}
In the next section we outline all the necessary results from operator theory that we need to prove the factorization theorem and to provide the new proof of the ratios formulas. This section requires very little background knowledge of operator theory and is nearly self-contained. The basic properties of Toeplitz operators are presented along with a proof of the Borodin-Okounkov-Case-Geronino (BOCG) identity.  

The next  section describes the connections with the  averages over the classical groups and determinants of finite sections of either Toeplitz or ``Toeplitz plus Hankel''.
This is followed by the proofs of the factorizations of the averages illustrated by the previous examples. Then Section 6 contains the new proofs of the ratio formulas using Fredholm determinant expansions. Finally, more examples can be found in the Appendix.

\section{Toeplitz  Preliminaries}
We begin with a function $\phi$ defined on the unit circle $\T$ with Fourier coefficients
$$\phi_{k} = \frac{1}{2\pi} \int_{0}^{2\pi} \phi(e^{i\theta})\,e^{-ik\theta}\,\,d\theta ,$$
$$ \phi(e^{i\theta}) = \sum_{-\infty}^{\infty} \phi_{k} \,e^{i\,k\,\theta} =  \sum_{-\infty}^{\infty} \phi_{k} \,z^{k},$$
($z = e^{i\theta} $) and consider the matrix 

$$T_n(\phi)\, = \, (\phi_{j-k})_{j,\,k\, = \,0,\, \cdots,\, n-1}.$$

The matrix $T_{n}(\phi)$ is referred to as the finite Toeplitz matrix with symbol $\phi.$

This matrix has the form 

\[ \left[\begin{array}{ccccc}\phi_{0} & \phi_{-1} & \phi_{-2} &  \cdots&  \phi_{-(n-1)} \\
\phi_{1} & \phi_{0} & \phi_{-1} &  \cdots  & \phi_{-(n-2)} \\
\phi_{2} & \phi_{1} & \phi_{0} & \cdots  & \phi_{-(n-3)} \\
\cdots & \cdots & \cdots & \cdots  & \cdots \\
\cdots & \cdots & \cdots&\cdots & \cdots \\
\phi_{n-1} & \phi_{n-2} & \phi_{n-3} & \cdots & \phi_{0}
\end{array}\right].
\]

%
%
%
%
%
%
%
%
 
Since the 1950s there has been considerable interest in the asymptotics of the determinants of Toeplitz matrices as $n \rightarrow \infty.$ The interest was motivated by the study of the two-dimensional Ising model correlations, but other important connections are related to many other areas of applied mathematics, such as entanglement computations in statistical physics or problems in random matrix theory, See for example \cite{IJK, IMM, JK} for some applications and \cite{BS99} and \cite{BS06} for more history and general facts about Toeplitz operators.

For a certain class of ``nice'' symbols the asymptotics are described by the 
Szeg\"o-Widom Limit Theorem. This theorem states that if the symbol
$\phi$ defined on the unit 
circle $\T$ has a sufficiently well-behaved logarithm then the determinant of the Toeplitz
matrix
$$T_n(\phi)=(\phi_{j-k})_{\,\,j,\,k\,=\,0,\, \cdots,\,n-1}$$
has the asymptotic behavior
$$  D_{n}(\phi) := \det T_n(\phi)  \, \sim G(\phi)^n\,E(\phi)\ \ \ {\rm as}\ n\to\iy.$$ 
 
Here are the constants:
$$G(\phi) = e^{(\log \phi)_0}$$
and
$$ E(\phi)  = \det \left( T(\phi)T(\phi^{-1}) \right)$$ 
where $\phi^{-1} = 1/\phi$ and
$$T(\phi) = (\phi_{j-k})\,\,\,\,\,\,\,\,0 \leq j, k < \iy $$ 

is the Toeplitz operator defined on $\ell^{2},$ the space of complex sequences $\{x_k\}_{k=0}^\infty,$ with the usual $2$-norm.

If $\phi$ is in $L^{\infty}(\T)$ then it is well known that $T(\phi)$ is a bounded operator on $\ell^{2}.$ More about the properties of Toeplitz operators will be described later. 
 
$E(\phi)$ is often referred to as Widom's Constant since Widom discovered this form of the constant when he extended the original Szeg\"o theorem to the block matrix case.

To make sense of Widom's constant, that is the term $\det \left( T(\phi)T(\phi^{-1}) \right)$,
we need to think about operators of the  form $I + K$ where $K$ is trace class. A compact operator $K$ is trace class if 
the eigenvalues $s_{i}$ of $(K^{*}K)^{1/2}$ satisfy $\sum_{i=0}^{\infty} s_{i} < \infty.$  When this is finite,  we define the trace norm of $K$ as
\[ ||K||_{1} = \sum_{i=0}^{\infty} s_{i}.\] 
This set of operators is a closed ideal in the set of all bounded operators. 

For trace class operators with discrete eigenvalues $\lambda_{i}$ one can show that  
\[ \sum_{i=0}^{\infty}|\lambda_{i}| < \infty \] and thus
\[ \det (I + K) = \prod_{i=0}^{\infty} (1+ \lambda_{i}) \]
is well defined.
Proofs of the above statements can be found in \cite{GK}.

We will later show that the operator $T(\phi)T(\phi^{-1})$ is of this form. 

To state the necessary conditions for the Szeg\"o-Widom theorem to hold, it is useful to consider a certain Banach algebra (see \cite{BS06}, section 1.10).

Let $B$ stand for the set of all function $\phi$ such that the Fourier coefficients satisfy
\[ \| \phi\|_{B} :=\sum_{k =-\iy} ^{\iy} |\phi_{k}| + 
 \Big(\sum_{k = -\iy}^{\iy} |k|\cdot|\phi_{k}|^{2}\Big)^{1/2} < \iy. \]
With this norm  and pointwise defined algebraic operations on $T$ the set $B$ becomes a Banach algebra of continuous functions on the unit circle.

The Szeg\"o-Widom Limit Theorem  holds provided that
$\phi\in B$ and the function $ \phi$ does  not vanish on $\T$ and has winding number zero. These conditions imply that both $\phi^{-1}$ and $\log \phi$ (any continuous choice of the logarithm will suffice) are defined and in $B.$

%
%

The most direct way to prove the Szeg\"{o}-Widom theorem is to prove an identity for the determinants, an identity called the Borodin-Okounkov-Case-Geronimo (BOCG) identity. We should point out here that while the identity easily yields the asymptotics, the asymptotics were known long before the identity.


In addition to the Toeplitz operator, we also define a Hankel operator
\begin{eqnarray*}
 T(\phi) & =&  (\phi_{j-k}), \,\,\,\,\,\,\,\,\,\,\,\,\,0 \leq j,k < \infty,\\
 H(\phi) &=&  (\phi_{j+k+1}), \,\,\,\,\,\,\,\,0 \leq j,k < \infty.
 \end{eqnarray*}
 
 It is useful to think of the Toepltiz operator in matrix form, but it is also worth pointing out that it can be defined as 
 \[ T(\phi) : \ell^{2} \rightarrow \ell^{2} \,\,\,\, \quad T(\phi)f = P(\phi f) \] 
 where $P$ is the Riesz projection onto $\ell^{2},$ 
 \[ P : \sum_{k=-\infty}^{\infty}\phi_{k}e^{ik\theta} \rightarrow \sum_{k = 0}^{\infty}\phi_{k}e^{ik\theta} ,\] and $f$ is the function identified by its Fourier coefficients $\{f_{n}\}_{n=0}^{\infty}.$ In other words, start with the sequence $\{f_{n}\}_{n = 0}^{\infty}$ in $\ell^{2}$  and then let 
 \[ f(\theta)  = \sum_{n=0}^{\infty} f_{n} e^{i n \theta} ;\] multiply pointwise by the function $\phi$ and then project by chopping off all negative Fourier coefficients of $\phi f$ to obtain a sequence back in $\ell^{2}.$

For $\phi, \psi \in L^{\infty}(\mathbb{T})$ the well-known identities
\begin{eqnarray}\label{talpro}
T(\phi \psi) & =& T(\phi)T(\psi) +H(\phi)H({\tilde \psi})\\
H(\phi \psi) & =& T(\phi)H(\psi) +H(\phi)T({\tilde \psi})
 \end{eqnarray}
 where $\tilde{\phi}( e^{i\,\theta}) = \phi (e^{-i\,\theta})$
are the basis of almost everything we do later on.  
They can be proved by simply analyzing the matrix entries of the operators and using convolution to find the Fourier coefficients of products. They are also found in \cite{BS06}, Proposition 2.14.

It follows from these identities that if  $\psi_{-}$ and $\psi_{+} $ have the property that all their Fourier coefficients vanish for $k > 0$ and $k < 0$, respectively, then 
\begin{eqnarray*}
T(\psi_{-} \phi \psi_{+} )& =& T(\psi_{-} )T(\phi) T( \psi_{+} ) ,\\
& &\\
H(\psi_{-} \phi \tilde{\psi}_{+} )& =& T(\psi_{-} )H(\phi) T( \psi_{+} ).
\end{eqnarray*} 

Now suppose that the conditions of the Szeg\"o-Widom theorem hold. Then there exists a well-defined logarithm of $\phi$ that is in $B$. We can then  write $\log \phi$ as a sum of two functions, one with Fourier coefficients that vanish for positive indices, and one with Fourier coefficients that vanish for negative indices. If we then exponentiate these two functions, the factors (and their inverses) have the same Fourier properties and are in $B$. Thus we see that
\[ \phi = \phi_{-}\phi_{+}, \,\,\, \phi^{-1} = \phi_{-}^{-1} \phi_{+}^{-1} .\]  This factorization is known as a Wiener-Hopf factorization. 

The BOCG identity can now be stated as follows.
\begin{theorem}\label{BOCG} Suppose $\phi$ satisfies the conditions of the Szeg\"o-Widom theorem. Then
\begin{eqnarray*}
D_{n}(\phi) &=& G(\phi)^n \,\,E(\phi) \cdot 
\det\left(I-H(z^{-n}\phi_{-}\phi_{+}^{-1})H(\tilde{\phi}_{-}^{-1} \tilde{\phi}_{+}z^{-n})\right).
\end{eqnarray*}
In the above $z = e^{i\theta}.$
\end{theorem}
\proof

We let 
$P_{n}$ be the projection defined on $\ell^{2}$ that maps $$\{x_k\}_{k=0}^\infty \,\,\,\mbox{to}\,\,\,\{x_0,\dots,x_{n-1},0,0, \dots\}$$ identifying a function with its Fourier coefficients, and we define $Q_{n} = I - P_{n}.$

It is easy to check that $$Q_{n} = T(z^{n})T(z^{-n}), \,\,\,I =  T(z^{-n})T(z^{n}).$$

In what follows and for the rest of the paper, determinants are always defined on the appropriate space. Thus for example, 
$\det P_{n} A P_{n}$ is the determinant defined on the image of $P_{n}.$ 

Next we make some simple observations:

1) If $\phi_{k} = 0$ for $k < 0$ or $k > 0, $ then the Toeplitz matrices are triangular and $D_{n}(\phi) = \det T_{n}(\phi) = (\phi_{0})^{n}.$

2) $T_{n}(\phi) = P_{n}T(\phi)P_{n}$ i.e. $T_{n}(\phi)$ is the upper left corner of the Toeplitz operator $T(\phi).$

The projection $P_{n}$ has the nice property that if
$U$ is an operator whose matrix representation has an upper triangular form, then
\[P_{n} U P_{n} = U P_{n}.\] 
And if $L$ is an operator whose matrix representation has a lower triangular form, then
\[ P_{n} L P_{n} =  P_{n} L. \]

 So if we had an operator of the form $LU,$ then
\[P_{n}LU P_{n} = P_{n}L P_{n}UP_{n}\]
and the corresponding determinants would be easy to compute.

What happens for Toeplitz operators is the opposite. 

Toeplitz operators generally do not factor as $L\,U$ but rather as $U\,L.$

To see this, we use a Wiener-Hopf factorization. We have 
\[\phi = \phi_{-} \,\phi_{+},\,\,\,\,\ \phi^{-1} = \phi_{-}^{-1}\phi_{+}^{-1}.\] 
Recall this means that  $(\phi_{+})_{k} = 0$ for $k < 0$ and thus
$T(\phi_{+})$ is lower triangular
and that
 $(\phi_{-})_{k} = 0$ for $k >0$ and $T(\phi_{-})$ is upper triangular.
We also know from the algebra property (\ref{talpro}) that  the factorization of $\phi$ shows that \[T(\phi) = T(\phi_{-}) \,T(\phi_{+}) .\]
One final observation, the algebra property (\ref{talpro}) implies that
\[ T(\phi_{+}^{-1})T( \phi_{+}) =   T(\phi_{+})T( \phi_{+}^{-1}) = I\]
and that
\[ T(\phi_{-}^{-1})T( \phi_{-}) =   T(\phi_{-})T( \phi_{-}^{-1}) = I.\]

We write
\begin{eqnarray*}P_{n} T(\phi) P_{n} &=& P_{n} T(\phi_{-}) \,T(\phi_{+}) P_{n}\\
&=& P_{n} T(\phi_{+}) T(\phi_{+}^{-1})T(\phi_{-}) \,T(\phi_{+}) T(\phi_{-}^{-1}) T(\phi_{-})P_{n} \\
  &=& P_{n} T(\phi_{+}) P_{n}T(\phi_{+}^{-1})T(\phi_{-}) \,T(\phi_{+}) T(\phi_{-}^{-1}) P_{n}T(\phi_{-})P_{n}.
  \end{eqnarray*}

Taking determinants, we have $$\det P_{n} T(\phi_{+}) P_{n} = ((\phi_{+})_{0})^{n}$$ and $$\det P_{n}T(\phi_{-}) P_{n} = ((\phi_{-})_{0})^{n}.$$
 So at this point we have
 \[ D_{n}(\phi) = G(\phi)^{n} \det  P_{n}T(\phi_{+}^{-1})T(\phi_{-}) \,T(\phi_{+}) T(\phi_{-}^{-1}) P_{n} .\]
   
The term 
\[ P_{n}T(\phi_{+}^{-1})T(\phi_{-}) \,T(\phi_{+}) T(\phi_{-}^{-1}) P_{n} \] 
is of the form 
\[ P_{n} (I + K) P_{n}\]
where $K$ is trace class. To see this we note that an operator is Hilbert-Schmidt if, in matrix form, its entries satisfy
\[ \sum_{i,j \geq 0} |a_{i,j}|^{2} < \infty ,\] and that it is well-known that the product of two Hilbert-Schmidt operators is trace class. For a good reference for this fact and general properties of trace class and Hilbert-Schmidt operators, please see \cite{GK}.

Now the operator 
$$ T(\phi_{+}^{-1})T(\phi_{-}) = T(\phi_{+}^{-1}\phi_{-}) - H(\phi_{+}^{-1}) H(\tilde{\phi_{-}}) $$ and the Banach algebra norm condition means that each of the above Hankel operators are Hilbert-Schmidt and thus their product is trace class. 
Thus
$$ T(\phi_{+}^{-1})T(\phi_{-}) = T(\phi_{+}^{-1}\phi_{-}) + K_{1} $$
where $K_{1}$ is trace class.
Applying the same argument to $T(\phi_{+}) T(\phi_{-}^{-1})$
we have
$$T(\phi_{+}) T(\phi_{-}^{-1}) = T(\phi_{+}\phi_{-}^{-1}) + K_{2}$$ where $K_{2}$ is trace class. 
Hence
$$T(\phi_{+}^{-1})T(\phi_{-})T(\phi_{+}) T(\phi_{-}^{-1})  =  (T(\phi_{+}^{-1}\phi_{-}) + K_{1})(T(\phi_{+}\phi_{-}^{-1}) + K_{2})$$
$$= T(\phi_{+}^{-1}\phi_{-})T(\phi_{+}\phi_{-}^{-1}) +  T(\phi_{+}^{-1} \phi_{-})K_{2} + K_{1}T(\phi_{+}\phi_{-}^{-1}) + K_{1}K_{2}.$$
This yields for a trace class $K_{3}$
$$ T(1) + K_{3} + T(\phi_{+}^{-1} \phi_{-})K_{2} + K_{1}T(\phi_{+}\phi_{-}^{-1}) + K_{1}K_{2}.$$
$$ = I +K$$ since $T(1) = I$ and trace class operators form an ideal. 

To obtain the next step we use the Jacobi identity 
\[ \det P_{n}A P_{n} = \det A\cdot\det (Q_{n}A^{-1} Q_{n})\]
in which $Q_{n}=I-P_{n}$ and $A$ is an invertible operator of the form identity plus trace class.

This holds because
\begin{eqnarray*}
\det P_{n}A P_{n} &=& \det (Q_{n}+P_{n}A P_{n})=\det (Q_{n}+A P_{n})\nn\\
&=& \det A\cdot \det (A^{-1} Q_{n}+P_{n})=
\det A\cdot\det (Q_{n}A^{-1} Q_{n}).\nn
\end{eqnarray*}

 We apply this to $A =  T(\phi_{+}^{-1})T(\phi_{-}) \,T(\phi_{+}) T(\phi_{-}^{-1}) $. It  follows that
\begin{eqnarray*}
\lefteqn{\det P_{n}T(\phi_{+}^{-1})T(\phi_{-}) \,T(\phi_{+}) T(\phi_{-}^{-1})P_{n}}
\hspace{6ex}\\[1ex]
&=&
\det T(\phi_{+}^{-1})T(\phi_{-}) \,T(\phi_{+}) T(\phi_{-}^{-1})\\&&\qquad \times   
\det Q_{n}T(\phi_{-})T(\phi_{+}^{-1}) \,T(\phi_{-}^{-1} )T(\phi_{+})Q_{n}.
\end{eqnarray*}
Using the fact that $\det AB = \det BA$ the first term on the right is easily seen to be $\det T(\phi)T(\phi^{-1}) = E(\phi).$ 
Thus we now have 
\[ D_{n}(\phi) = G(\phi)^{n} E(\phi) \det Q_{n}T(\phi_{-})T(\phi_{+}^{-1}) \,T(\phi_{-}^{-1} )T(\phi_{+})Q_{n}.\]
Now 
\begin{flalign*}
 \det Q_{n}A^{-1}Q_{n} &  =\det Q_{n}T(\phi_{-})T(\phi_{+}^{-1}) \,T(\phi_{-}^{-1} )T(\phi_{+})Q_{n}&  \\ 
 &  = \det (P_{n }+Q_{n}T(\phi_{-}\phi_{+}^{-1}) \,T(\phi_{-}^{-1} \phi_{+})Q_{n})& \\
 & =  \det( P_{n} + Q_{n}(I - H(\phi_{-}\phi_{+}^{-1} )H(\tilde{\phi}_{-}^{-1} \tilde{\phi}_{+}))Q_{n} )&\\
 & = \det( I - Q_{n}H(\phi_{-}\phi_{+}^{-1} )H(\tilde{\phi}_{-}^{-1} \tilde{\phi}_{+})Q_{n} )&\\
 & = \det( I - T(z^{n})T(z^{-n})H(\phi_{-}\phi_{+}^{-1} )H(\tilde{\phi}_{-}^{-1} \tilde{\phi}_{+})T(z^{n})T(z^{-n}) )&\\
 & = \det( I - T(z^{-n}) H(\phi_{-}\phi_{+}^{-1} )H(\tilde{\phi}_{-}^{-1} \tilde{\phi}_{+})T(z^{n} ))&\\
 & = \det( I - H(z^{-n}\phi_{-}\phi_{+}^{-1} )H(\tilde{\phi}_{-}^{-1} \tilde{\phi}_{+}z^{-n}) )
\end{flalign*}
and the identity is proved.
\endproof
For other proofs (including two original ones) we refer the reader to \cite{BW, BO, Bot, GC}.

From this BOCG identity we have an instant proof of the Szeg\"{o}-Widom theorem. Notice that $Q_{n}f$ tends to zero for any fixed $f$ as $n \rightarrow \infty$ and this allows one to say that 
$$Q_{n}\,H(\phi_{-}\phi_{+}^{-1})H(\tilde{\phi}_{-}^{-1} \tilde{\phi}_{+})\,Q_{n}$$ tends to zero in the trace norm. Thus the determinant
$$ \det ( I  - Q_{n}\,H(\phi_{-}\phi_{+}^{-1})H(\tilde{\phi}_{-}^{-1} \tilde{\phi}_{+})\,Q_{n}) $$ tends to one and 
$$  D_{n}(\phi) = \det T_n(\phi)  \, \sim G(\phi)^n\,E(\phi)\ \ \ {\rm as}\ n\to\iy.$$
For more details see \cite{BW}.

%
%
%
%
%
 
The constant $E(\phi)$ has a nice concrete description. 

If we have a Wiener-Hopf factorization for $\phi = \phi_{-}\phi_{+}$, then 
$$ T(\phi)T(\phi^{-1}) = T(\phi_{-})T(\phi_{+})T^{-1}(\phi_{-})T^{-1}(\phi_{+})$$ and this is of the form
$$e^{\,A}e^{\,B}e^{-A}e^{-B}$$ 
where $$A = T(\log(\phi_{-})), \,\,\,\,B = T(\log (\phi_{+})).$$

From this we can use a formula for determinants of multiplicative commutators of this form,
$$\det \,(e^{A}e^{B}e^{-A}e^{-B}) = \exp\,\,(\tr (AB - BA)).$$ 

It is tempting to think that the above formula is trivial, that is, that the left and right sides should be $1$ since the trace of $AB$ should be the trace of $BA$ and that the operators on the left side should be inverses of each other. But this is not the case. For our $A$ and $B,$ $AB$ is not trace class (indeed not even compact), but from our Banach algebra conditions what is true is that $AB - BA$ is trace class. Moreover, the operators $e^{A}$ and $e^{B}$ do not commute. In this situation the formula above yields a non-trivial answer. This was first used in this context in \cite{Wi76}.

Now $AB - BA$ is the same as
\[ T(\log(\phi_{-}))T(\log (\phi_{+})) - T(\log (\phi_{+}))T(\log(\phi_{-})) \]
\[ = T(\log(\phi_{-})\log (\phi_{+})) - T(\log (\phi_{+}))T(\log(\phi_{-}))\] which by \ref{talpro}
is \[ H(\log (\phi_{+}))H(\widetilde{\log (\phi_{-})}) .\] After computing the trace and exponentiating we arrive at 
$$\det T(\phi)T(\phi^{-1}) =  \exp\left(\sum_{k=1}^{\iy}k\,(\log\phi)_k\,(\log\phi)_{-k}\right).$$
 
Now we turn to a different class of structured operators.

\section{Toeplitz plus Hankel operators}

Many other classes of structured matrices arise in applications other than finite Toeplitz matrices. For example, one may want to find asymptotics for the determinants of matrices generated by the Fourier coefficients of $\phi$ with $j,k$ entry
\[ \phi_{j-k} + \phi_{j+k+1} .\]

The purpose of this section is to state identities similar to the BOCG type proved in the previous section. The proofs of these are similar (but not exactly the same) as the one done in the previous section and the proofs and additional details can be found in \cite{BE4}.

Let $\cS$ stand for a unital Banach algebra of functions on the unit circle 
continuously embedded into $L^\iy(\T)$ with norm $||\cdot ||_{S}$ and with the property that 
$\phi \in \cS$ implies that $\tilde{\phi}\in\cS$ and $P \phi\in \cS$. Moreover, define
\bq
\cS_- &=& \Big\{ \phi\in \cS\;:\; \phi_n=0 \mbox{ for  all } n>0\Big\},\nn\\
\cS_{+}&=&\Big\{ \phi\in \cS\;:\; \phi_n=0 \mbox{ for  all } n<0\Big\},\nn\\
\cS_0 &=& \Big\{ \phi\in \cS\;:\; \phi=\tilde{\phi} \Big\}.\nn
\eq

In the above $\tilde \phi (e^{i\theta}) = \phi(e^{-i\theta})$
and $P$ is the Riesz projection. 
 
Assume that $M: \phi \in L^\iy \mapsto M(\phi)\in \cL(\ell^2)$ (the set of bounded operators on $\ell^{2}$) is a continuous linear map such that:
\begin{enumerate}
\item[(a)]
If $\phi \in \cS$, then $M(\phi)-T(\phi)$ is trace class
and there exists some constant $C$ such that
$$ 
\|M(\phi)-T(\phi)\|_{1} \le C\, \|\phi\|_{\cS}.
$$
\item[(b)]
If $\psi \in \cS_-$, $\phi \in \cS$, $ \gamma \in \cS_0$, then
$$
M(\psi \phi \gamma)=T(\psi)M(\phi) M(\gamma).
$$
\item[(c)]
$M(1)=I$.
\end{enumerate}
Then we say $M$ and $\cS$ are compatible pairs.
 
All of the following can be realized as compatible pairs.


\begin{enumerate}
\item[(I)]
$M(\phi) =T(\phi)+H(\phi)$,
\item[(II)]
$M(\phi) =T(\phi)-H(\phi)$,
\item[(III)]
$M(\phi)=T(\phi)-H(z^{-1} \phi)$, 
\item[(IV)]
$M(\phi)=T(\phi)+H( z \,\phi)Q_{1}$.
\end{enumerate}

The matrix representations of the operators are of the form
\[ (\phi_{j-k} \pm \phi_{j+k-\kappa+1}) \]
with $\kappa = 0, 1, -1.$

A Banach algebra that satisfies the compatible pair conditions for all four of the above cases is the set of all $\phi$ such that 
\[ |\phi_{0}| + ||H(\phi)||_{1}+ ||H(\tilde{\phi})||_{1} < \infty.\]

We are interested in the determinants (where the matrices or operators are always thought of as acting on the image of the projection of the appropriate space) of
\[ P_{n}\, M(\phi) \,P_{n},\]
or for our examples, 
determinants of finite Toeplitz plus finite Hankel matrices.

%

 
%
%
%
%

As mentioned earlier the analog of the BOCG for the determinants of $P_{n}M(\phi)P_{n}$ was proved in \cite{BE4}. Here is the statement of the identity for even functions ($\phi = \tilde{\phi}$), but can be made more general. Note that in the theorem below $E(\phi)$ always denotes the constant in the identity and is different depending on the compatible pair.

\begin{theorem}\label{BEf}
Let $M$ and $\cS$ be a compatible pair, and let $\beta_+\in \cS_+$. Let $\phi=\phi_+\tilde{\phi}_+=\exp(\beta)$ with  $\phi_+=\exp (\beta_+)$, $\beta=\beta_++\tilde{\beta}_+$.
Then
\[ \det P_{n} M(\phi) P_{n}  = G(\phi)^{n} E(\phi) \det (I + Q_{n} \,K \,Q_{n} ),\]
where
\[ E(\phi) =  \exp\Big( \tr(M( \beta)-T(\beta)) + \frac{1}{2}\tr\, H(\beta)^2\Big),\]
and $K = M(\phi_{+}\iv) T(\phi_{+}) - I.$
\end{theorem}

For case (I),\,\,\,\, $(\phi_{j - k}  + \phi_{j+k + 1}),$
$$K = M(\phi_{+}\iv) T(\phi_{+}) - I = H(\phi_{+}^{-1})T(\phi_{+}) = H(\tilde \phi_{+}/\phi_{+}).$$
For case (II),\,\,\,\,$(\phi_{j - k}  - \phi_{j+k + 1}),$
$$K = M(\phi_{+}\iv) T(\phi_{+}) - I = -H(\phi_{+}^{-1})T(\phi_{+}) = -H(\tilde \phi_{+}/\phi_{+}).$$
For case (III),\,\,\,\, $(\phi_{j - k} - \phi_{j+k + 2}),$
$$K = M(\phi_{{+}}\iv) T(\phi_{+}) - I,$$
$$M(\phi)=T(\phi)-H(z^{-1} \phi),$$ we have 
\begin{eqnarray}\nonumber
 K &= &(T(\phi_{+}^{-1}) - H(z^{-1}\phi_{+}^{-1} ))T(\phi_{+}) - I\\ \nonumber
& = &T(\phi_{+}^{-1})T(\phi_{+})  - H(z^{-1}\phi_{+}^{-1} ))T(\phi_{+}) - I\\ \nonumber
& = & - H(z^{-1}\phi_{+}^{-1} ))T(\phi_{+})\\ \nonumber
& = & - H(z^{-1}\phi_{+}^{-1} \tilde{\phi_{+}}).
\end{eqnarray}

For case (IV), 
$$K = M(\phi_{+}\iv) T(\phi_{+}) - I,$$
$$M(\phi)=T(\phi) + H(z \phi)Q_{1},$$ we have 
\begin{eqnarray}\nonumber
 K &= &(T(\phi_{+}^{-1}) + H(z\,\phi_{+}^{-1} )Q_{1})T(\phi_{+}) - I\\ \nonumber
& = &T(\phi_{+}^{-1})T(\phi_{+})  + H(z\, \phi_{+}^{-1} )Q_{1}T(\phi_{+}) - I\\ \nonumber
& = &  H(z\, \phi_{+}^{-1} )Q_{1}T(\phi_{+}).
\end{eqnarray}
Note for case (IV), 
$$M(\phi)\, R = (\phi_{j-k} + \phi_{j+k})$$ where $R = 2P_{1} + Q_{1},$
that is, a diagonal matrix with $2$ in the upper-left entry and ones otherwise.

Note also that 
$$ \det P_{n}\, M(\phi) \, R\, P_{n} = 2 \det P_{n}\,M(\phi) \, P_{n}. $$
The need for formulas for the above determinants and their asymptotics arise in various applications. One example is the connection to averages over certain groups or cosets of groups. We will now state the connections. The measure in all cases is normalized Haar measure. All of the results are found in Section 2 of \cite{BR}.
\begin{proposition} ${}$
Suppose $ \phi$ and $g$  satisfy the conditions of Theorem \ref{BEf}. Let $N = 2m.$  Let $N = 2m.$ Then
\begin{enumerate}

\item Let $I_{U(N)}(\phi)$ denote the average over the unitary group $U(N)$  of $\det \phi(U).$ Then
$$I_{U(N)}(\phi) = D_{N}(\phi).$$
\item Let $N = 2m, \phi = g \tilde g$. Then 
$$I_{O^{-}(N+1)}(g) = g(-1)\det P_{m}(\phi_{j-k} + \phi_{j+k+1})P_{m}$$ averaging $\det g(U)$ over the coset of orthogonal matrices with determinant $-1.$
\item Let $N=2m, \phi =g \tilde g$.Then $$I_{O^{+}(N+1)}(g) = g(1)\det P_{m}(\phi_{j-k} - \phi_{j+k+1})P_{m}$$ averaging $\det g(U)$ over the subgroup of matrices with determinant $1.$
\item Let $N = 2m, \phi = g \tilde g$ Then
$$I_{USp(N)}(g) = \det P_{m}( \phi_{j-k} - \phi_{j+k+2}) P_{m} $$ averaging  $\det g(U)$ over the symplectic group.
\item Let $N = 2m, \phi =g \tilde g.$ Then $$I_{O^{+}(N)}(g) =  \frac{1}{2}\det P_{m}(\phi_{j-k} + \phi_{j+k}) P_{m} $$ averaging $ \det g(U)$ over the subgroup of matrices with determinant $1.$
\item Let $N = 2m, \phi = g\tilde g.$ Then 
$$I_{O^{-}(N)}(g) =  g(1)g(-1)\det P_{m-1}(\phi_{j-k} - \phi_{j+k+2}) P_{m-1} $$ averaging $ \det g(U)$ over the coset of matrices with determinant $-1.$

\end{enumerate}
Note in the formulas above  $g(\pm 1)$  means substitute $\pm 1$  for $e^{i\theta}$ in the formula for $g$. So for
example, if
$g = 1 - a e^{i\theta}$  then $g(1) = 1 - a$ and $g(-1)=1+a$.

\end{proposition}

In the following sections we shall see how the above formulas and the determinant identities of this section arise in number theory calculations. 

\section{Factoring determinants}

The example given in the introduction showed a factorization of the averages over the unitary groups for the $\det \phi$ when $\phi$ was a certain rational function. It turns out that is always the case for even functions. 

The following theorem describes the factorization of the Toeplitz determinant when $N$ is even.
\begin{theorem}
Suppose $ \phi = g \tilde g$ satisfies the conditions of Theorem \ref{BEf}. Let $N = 2m.$
Then 
\[ I_{U(N)}(g\tilde{g}) = (g(1)g(-1))^{-1}I_{O^{+}(N+1)}(g) I_{O^{-}(N+1)}(g),\] or in the operator language
$ D_{N}(\phi) $ equals \[ \det P_{m}(T(\phi) + H(\phi)) P_{m} \times \det P_{m}( T(\phi) - H(\phi)) P_{m}.\]

\end{theorem}
\proof
The proof of this follows directly from the formulas of Theorems 3.1 and \ref{BEf}. It is easy to see that both the geometric mean part and the constant term part of the determinants on the right-hand side of the above multiply to give the correct term for the Toeplitz or Unitary case. For the Fredholm determinant part of the identity we have
$$\det ( I - Q_{N} \,K\, Q_{N}) $$ where 
$$K = H( \tilde{\phi_{+}}/\phi_{+} )^{2}.$$
Now 
$$Q_{N}H( \tilde{\phi_{+}}/\phi_{+} )^{2} Q_{N} = T(z^{2m})T(z^{-2m})H( \tilde{\phi_{+}}/\phi_{+} )^{2}T(z^{2m})T(z^{-2m}).$$  
This is the same as
$$T(z^{2m})T(z^{-m})H( z^{-m}\tilde{\phi_{+}}/\phi_{+} )H(\tilde{\phi_{+}}/\phi_{+} )T(z^{m})T(z^{m})T(z^{-2m})$$  or
\begin{eqnarray*} &&T(z^{2m})T(z^{-m})H( z^{-m}\tilde{\phi_{+}}/\phi_{+} )T(z^{-m})\\ &&\qquad \times T(z^{m})H(\tilde{\phi_{+}}/\phi_{+}) T(z^{m})T(z^{m})T(z^{-2m}),\end{eqnarray*}
which equals 
$$T(z^{2m})T(z^{-m})H(\tilde{\phi_{+}}/\phi_{+}) )T(z^{m})T(z^{-m})H(\tilde{\phi_{+}}/\phi_{+}) T(z^{m})T(z^{-2m}).$$ 
Hence the determinant of $ I - Q_{N} \,K_{u}\, Q_{N}$ is the same as
\[ \det( I - T(z^{m}) Q_{m} H(\tilde{\phi_{+}}/\phi_{+})Q_{m} H(\tilde{\phi_{+}}/\phi_{+}) Q_{m} T(z^{-m})) \]
or
\begin{eqnarray*} && \det (I -  Q_{m} H(\tilde{\phi_{+}}/\phi_{+}) Q_{m} H(\tilde{\phi_{+}}/\phi_{+}) Q_{m} ) \\ && \qquad 
= \det (I - Q_{m} H(\tilde{\phi_{+}}/\phi_{+}) Q_{m}) \det (I + Q_{m} H(\tilde{\phi_{+}}/\phi_{+}) Q_{m})\end{eqnarray*}
and the result follows.
\endproof
Next we have the case when $N$ is odd.
\begin{theorem} 
Suppose $\phi = g \tilde g$ satisfies the conditions of Theorem \ref{BEf}. Let $N -1 = 2m,$
\[ I_{U(N)}(g\tilde{g}) = I_{USp(N-1)}(g) I_{O^{+}(N+1)}(g),\] or in the operator language
$ D_{N}(\phi) $ equals \[\det P_{m}(T(\phi)- H(z^{-1}\phi)) P_{m} \times \frac{1}{2}\det P_{m +1}( T(\phi) + H(z\phi)Q_{1})R\, P_{m+1}\]
\[ = \det P_{m}(T(\phi)- H(z^{-1}\phi)) P_{m} \times \det P_{m +1}( T(\phi) + H(z\phi)Q_{1})\, P_{m+1}.\]

\end{theorem}
\proof
As before, each formula has a $G(a)$ part, a constant part, and a Fredholm determinant part. 

For $I_{U(N)},$  the exponent of $G(a)$ is $N = 2m+1.$ 
For $I_{USp(N-1)}$ the exponent is $m$ and for $ I_{O^{+}(N+1)},$  $m+1.$

We write $\phi  =  \phi_{+}\tilde{\phi_{+}} = \exp \beta = \exp (\beta_{+} + \tilde{\beta_{+}}).$ For  $I_{U(N)},$ we have constant term $E(\phi) = \exp\,( \tr \,H(\beta)^{2})$. 

For $I_{USp(N-1)}$ the constant term is $$\exp \tr (-H(z^{-1}\beta) + \frac{1}{2} H(\beta)^{2}).$$
For $ I_{O^{+}(N+1)}$ it is $$\frac{1}{2} \times 2 \exp \tr (H(z \,\beta)Q_{1} + \frac{1}{2} H(\beta)^{2}).$$
But then the product of the last two is $$ \exp(- \sum_{j= 1}^{\infty} \beta_{2j} + \sum_{j = 1}^{\infty} \beta_{2j} + \tr \,H(\beta)^{2} ).$$

Now we compare the Fredholm determinants. 
For the unitary case,
we have $$\det ( I - Q_{N} \,K_{u}\, Q_{N}) $$ where 
$$K_{u} = (H( \tilde{\phi_{+}}/\phi_{+}) )^{2}.$$
For the symplectic case we have
$$ \det ( I - Q_{m} \, K_{s}\,Q_{m} )$$ with $$K_{s} = - H( z^{-1} \tilde{\phi_{+}}/\phi_{+}).$$
And for the orthogonal case we have
$$ \det ( I - Q_{m+1}\, K_{o} \,Q_{m+1} )$$ with $$K_{o} = H(z /\phi_{+}) Q_{1} T(\phi_{+}).$$

Thus we consider
$$ \det ( I  + Q_{m} \, H( z^{-1} \tilde{\phi_{+}}/\phi_{+})\,Q_{m} ) \det ( I - Q_{m+1}\, H(z /\phi_{+}) Q_{1} T(\phi_{+})\,Q_{m+1} ).$$
The determinant 
$$ \det ( I  + Q_{m} \, H( z^{-1} \tilde{\phi_{+}}/\phi_{+})\,Q_{m} )$$ is the same as
$$\det ( I + T (z^{m})T(z^{-m}) T(z^{-1}) H (\tilde{\phi_{+}}/\phi_{+})T(z^{m}) T(z^{-m}))$$ or
$$\det ( I + T (z^{m+1})T(z^{-m-1})  H (\tilde{\phi_{+}}/\phi_{+}) T(z^{m}) T(z^{-m-1}))$$ 
$$ \det ( I + Q_{m+1}  H (\tilde{\phi_{+}}/\phi_{+}) T(z^{m}) T(z^{-m-1})).$$

The determinant 
$$\det ( I - Q_{m+1}\,H(z /\phi_{+}) Q_{1} T(\phi_{+}) \,Q_{m+1} )$$ is the same as
$$ \det ( I -  Q_{m+1}H(z \,/\phi_{+}) T(z) T(z^{-1}) T(\phi_{+})T(z^{m+1})T(z^{-m-1}))$$ or
$$ \det ( I -  Q_{m+1}H( \tilde{\phi_{+}}/\phi_{+})  T(z^{m})T(z^{-m-1})).$$

Thus for the product we have $\det ( I - L)$ where $L$ is
$$  Q_{m+1} H(\tilde{\phi_{+}}/\phi_{+}) T(z^{m}) T(z^{-m-1})Q_{m+1}  H(\tilde{\phi_{+}}/\phi_{+}) T(z^{m}) T(z^{-m-1})$$
$$ Q_{m+1}H( \tilde{\phi_{+}}/\phi_{+}) T(z^{m})T(z^{-m-1})H( \tilde{\phi_{+}}/\phi_{+}) T(z^{m})T(z^{-m-1})$$
$$ T(z^{-m})Q_{2m+1}H( \tilde{\phi_{+}}/\phi_{+})H( \tilde{\phi_{+}}/\phi_{+}) T(z^{2m+1})T(z^{-m-1}),$$
and
$$\det ( I - T(z^{-m})Q_{2m+1}H( \tilde{\phi_{+}}/\phi_{+}) H( \tilde{\phi_{+}}/\phi_{+})  T(z^{2m+1})T(z^{-m-1}))   $$
$$= \det ( I - Q_{2m+1}H( \tilde{\phi_{+}}/\phi_{+}) H( \tilde{\phi_{+}}/\phi_{+}) Q_{2m+1}) $$
$$ = \det ( I - Q_{N}\, K_{u}\, Q_{N} ) .$$

Thus we have shown that when $N$ is odd, $D_{N}(\phi)$ is the same as
\[ \det P_{m}(T(\phi)- H(z^{-1}\phi)) P_{m} \times \det P_{m +1}( T(\phi) + H(z\phi)Q_{1})\, P_{m+1}\]
which is our desired result. 
\endproof
 
A factorization similar to the one found in Theorem 5.1 can be done for other classes of operators including
Wiener-Hopf and Wiener-Hopf plus/minus Hankel operators. Determinants of these operators are related to other ensembles of random matrices.

  \section{New proof of ratios formula}
   
  We now use the identities described in the previous sections, along with an expansion for a Fredholm determinant to prove the ratios formulas mentioned in the introduction.  
  In this section we use the following notation.
  
 For sets of complex numbers $A$ and $B$ let
$$Z(A,B)=\prod_{a\in A\atop b\in B}(1-ab)^{-1}.$$
Further,  let
$$Z(A,B;C,D)=\frac{Z(A,B)Z(C,D)}{Z(A,D)Z(B,C)}.$$

Also, define
 $$Z_S(A;C)=\frac{Z_S(A)Z_O(C)}{Z(A,C)}$$
 where, if $A=\{a_1,\dots,a_k\}$ then 
 $$Z_S(A)=\prod_{1\le i\le j\le k}(1-a_i a_j)^{-1}$$
 and 
 $$Z_O(A)=\prod_{1\le i< j\le k}(1-a_i a_j)^{-1}$$
 
For the set $A$ with $U \subset A,$ the set $A-U+U^{-1}$ is the set $A$ with the elements of $U$ removed and with  the inverses of the elements of $U$ put back in. Finally, $U^{N} = \prod_{a_{i} \in U}a_{i}^{N}.$
 
 \begin{theorem}\label{t1} Let \[g = \prod_{i=1}^{k}\frac{(1 - a_{i}e^{i \theta})}{(1 - c_{i}e^{i \theta})},\] and let $A = \{a_{1}, a_{2}, \dots, a_{k}\}$. Suppose that $|c_{i}| <1$ for all $i.$ Then 

$$I_{Sp(2m)}(g) = \sum_{U\subset A}U^NZ_S(A-U+U^{-1};C).$$
 
 \end{theorem}
  
\proof
We first note that the function $g$ satisfies the conditions of Theorem \ref{BEf} providing that $|a_{i}| <1$ for all $i$. We make that assumption in our computation below, but then by analytic continuation the result holds for all $a_{i}.$

 We use the identity in Theorem \ref{BEf} to compute the determinant corresponding to the symplectic case exactly.
 Our symbol is \[ \phi(\theta)= 
g \tilde g = \prod_{i=1}^{k}\frac{(1 - a_{i}e^{i \theta})(1 - a_{i}e^{-i \theta})}{(1 - c_{i}e^{i \theta})(1 - c_{i}e^{-i \theta})}.\]
 
We see for this symbol the geometric mean $G(\phi) = 1.$ To compute $E(\phi)$ is straightforward and is seen to be
\[ \frac{\prod_{i,j} (1 - a_{i}c_{j})}{\prod_{1\leq i\leq  j \leq k}( 1 - a_{i}a_{j})\prod_{1 \leq i < j \leq k}(1 - c_{i}c_{j})}.  \]

Now for the Fredholm determinant part of the identity, we use the standard Fredholm expansion. Recall that if we have a discrete kernel $K(i,j),$ for the operator $K,$ then the determinant of $I - K$ is given by a sum
\[1 +\sum_{l = 1 }^{\infty} \frac{(-1)^{l}}{l!} \left(\sum_{i_{1}, i_{2}, \cdots, i_{l}} \det K(i_{j}, i_{k})_{1 \leq j,k \leq l}\right). \]
The first term here is simply the trace of the operator $K$. The second term is a two-by-two determinant, and so on. 

In our case we need to apply this to the Hankel operator $$K_s=Q_{m}H(z^{-1}\tilde{\phi_{+}}/\phi_{+})Q_{m}.$$
 So the first thing we need are the positive Fourier coefficients of the function, $\tilde{\phi_{+}}/\phi_{+}.$ These are the Fourier coefficients of 
\[ \prod_{i=1}^{k}\frac{(1 - a_{i}e^{-i \theta})(1 - c_{i}e^{i \theta})}{(1 - a_{i}e^{i \theta})(1 - c_{i}e^{-i \theta})}.\] The Fourier $l$th coefficient is given by 
\[ \sum_{j = 1}^{k}\al_{j}a_{j}^{l-1}, \,\,\,\,\,\, \al_{j} = \prod_{i \neq j}\frac{1}{(a_{j} - a_{i})}\prod_{i}\frac{(1 - a_{i}a_{j})(a_{j}- c_{i})}{(1 -c_{i}a_{j})}.\]
This is found by writing the usual integral for the Fourier coefficients, then letting $z = e^{i\theta}$ on the contour $|z| = 1$ and computing the residues at infinity and the points $z = 1/a_{i}.$
Notice first that any term in the Fredholm expansion with index greater than $k$ must vanish. This is easy to see since if we expand the determinants in columns one would need to pick two columns with the same $a_{j}$'s and thus the determinant vanishes. 

Notice also the effect of the $Q_{m}$ operators is to simply cut off the first $m$ rows and columns of the Hankel matrix. 
Now we will compare the terms of the sum with the terms in the ratio expansion. Letting $U$ be the empty set, the zero'th  term in the ratio expansion is $Z_{S}(A,C)$ which is exactly the $E(\phi)$ term
\[ \frac{\prod_{i,j} (1 - a_{i}c_{j})}{\prod_{1\leq i\leq  j \leq k}( 1 - a_{i}a_{j})\prod_{1 \leq i < j \leq k}(1 - c_{i}c_{j})}.  \]

Now suppose that the set $U$ has cardinality one. Let $U = \{a_{1}\}.$ Consider the first term in the expansion of our determinant with $K = Q_{m}H(z^{-1}\tilde{\phi_{+}}/\phi_{+})Q_{m}.$ This is the trace of the operator $K.$ The operator $K$ is a sum of rank one operators (one for each $a_{j}$). One of these when $j = 1$ has trace 
\[ \alpha_{1}a_{1}^{2m+1}/(1 - a_{1}^{2}) = a_{1}^{2m+1}/(1-a_{1}^{2})\prod_{i \neq 1}\frac{1}{(a_{1} - a_{i})}\prod_{i}\frac{(1 - a_{i}a_{1})(a_{1}- c_{i})}{(1 -c_{i}a_{1})} .\]

After multiplying through by the $E(\phi)$ term, we have
\[-a_{1}^{2m}\frac{\prod_{i \neq 1,j} (1 - a_{i}c_{j})\prod_{i}(1- c_{i}/a_{1})}{(1-1/a_{1}^{2})\prod_{\substack{j \neq 1 \\1\leq i\leq  j \leq k}}( 1 - a_{i}a_{j})\prod_{1 \leq i < j \leq k}(1 - c_{i}c_{j})\prod_{i \neq 1}(1 - a_{i}/a_{1})} .\]

After adjusting for the minus sign, this is exactly the answer given by the ratios formula for $Sp(2m).$ And thus the trace term in the Fredholm expansion yields the sum of such terms, one for each set of cardinality one. 

Now for the general term in the Fredholm expansion, we let $n = l$ and suppose that $l \leq k.$ We need to compute the term 

\[ \sum_{i_{1}, i_{2}, \cdots, i_{l}} \det K(i_{j}, i_{h})_{1 \leq j,h \leq l}\,.\]

To do this, first recall that each column of $K(i_{j}, i_{h})_{1 \leq j,h \leq l}$ is a sum of the terms 
$\al_{j}a_{j}^{i_{j} +i_{h} + 1}.$ We can expand the determinant by the columns. Thus when two of the $a_{j}$'s are the same the determinant vanishes since each is a multiple of the other. It follows that we only get a contribution when the $a_{j}$'s are distinct. 

Now suppose that  we wish to compute the determinant for the set of $a_{j}$'s with $ j = 1, 2, \dots, l.$ First, we factor out the product 
\[ \prod_{j} \al_{j} = \prod_{j}\left(\prod_{i \neq j}\frac{1}{(a_{j} - a_{i})}\prod_{i}\frac{(1 - a_{i}a_{j})(a_{j}- c_{i})}{(1 -c_{i}a_{j})}\right)\]
or
\[\prod_{\substack{i\neq j\\ 1\leq j \leq l\\ 1 \leq i \leq k}}\frac{1}{(a_{j} - a_{i})}\prod_{\substack{1 \leq j \leq l\\1 \leq i\leq k}}\frac{(1 - a_{i}a_{j})(a_{j}- c_{i})}{(1 -c_{i}a_{j})}.\]
What is left is columns with entries of the form $a_{j}^{i_{g} +i_{h} + 1}.$ We can also factor out the product of the $a_{j}$'s to have columns that are simply of the form $a_{j}^{i_{g} +i_{h}}.$ We next expand this using the permutation definition of the determinants for some fixed choice of the columns. To begin we assume the $a_{j}$'s are in ascending order in the columns. 
($\sigma$ denotes a permutation of the set $\{1, 2, \cdots, l\}.)$ We have then
\[ \sum_{\sigma}sgn(\sigma)a_{1}^{i_{\sigma(1)} + i_{1}}a_{2}^{i_{\sigma(2) }+ i_{2}}\cdots  a_{l}^{i_{\sigma(l)} +i_{l}}.\]
Now pair each $a_{j}^{i_{j}}$ with the corresponding $a_{j'}$ which also has $i_{j}$ as its exponent and then sum. 
The result is 
\[ \sum_{\sigma}sgn(\sigma)(a_{1}a_{2}\cdots a_{l})^{2m}/((1 - a_{1}a_{\sigma(1)})(1 - a_{2}a_{\sigma(2)}) \dots (1 - a_{l}a_{\sigma(l)})).\]
Now note that this answer is independent of the order of the columns, so we could have started with any permutation of the 
$a_{j}$'s. There are $l!$ such permutations and each one yields the same answer. 

If we factor out the term $(a_{1}a_{2}\cdots a_{l})^{2m}$ we see this itself is another determinant of the matrix with entries 
$$ \frac{1}{1 - a_{i}a_{j}}\,\,\,\,\, 1\leq i, j, \leq l,$$ and this a Cauchy determinant after factoring out the $a_{j}$'s from each column. 

The result is that we have 
\[ (a_{1}a_{2}\cdots a_{l})^{2m} \frac{\prod_{i = 2}^{l}\prod_{j = 1}^{i-1}( 1/a_{i} - 1/a_{j})( a_{i} - a_{j})}{\prod_{i = 1}^{l}\prod _{j = 1}^{l}( 1/a_{j} - a_{i})}.\]

We now just need to collect our answers and remember that we must multiply through by $E(\phi)$.
Here is what we have so far:

\begin{enumerate}
\item 
$$E(\phi) =  \frac{\prod_{i,j} (1 - a_{i}c_{j})}{\prod_{1\leq i\leq  j \leq k}( 1 - a_{i,}a_{j})\prod_{1 \leq i < j \leq k}(1 - c_{i}c_{j})}.$$
\item
The product of the $\al_{j}$ terms is 
$$ \prod \al_{j} = \prod_{\substack{i\neq j\\ 1\leq j \leq l\\ 1 \leq i \leq k}}\frac{1}{(a_{j} - a_{i})}\prod_{\substack{1 \leq j \leq l\\1 \leq i\leq k}}\frac{(1 - a_{i}a_{j})(a_{j}- c_{i})}{(1 -c_{i}a_{j})}.$$
\item
The answer from the determinant part
\[ (a_{1}a_{2}\cdots a_{l})^{2m} \frac{\prod_{i = 2}^{l}\prod_{j = 1}^{i-1}( 1/a_{i} - 1/a_{j})( a_{i} - a_{j})}{\prod_{i = 1}^{l}\prod _{j = 1}^{l}( 1/a_{j} - a_{i})}.\]
\end{enumerate}
Note that we factored out the product of the $a_{j}$'s from both the numerator and the denominator in the matrix columns for the computation in (3). 

We now simplify the product of the above by considering three types of terms.
\begin{enumerate}
\item
The terms with only $c_{j}$'s. This comes only from the $E(\phi)$ term and is 
$$\frac{1}{\prod_{1 \leq i < j \leq k}(1 - c_{i}c_{j})}.$$
\item 
The terms with $c_{j}$'s and $a_{j}$'s. This comes from (1) and (2) above and is
$$\prod_{i,j} (1 - a_{i}c_{j})\prod_{\substack{1 \leq j \leq l\\1 \leq i\leq k}}\frac{(a_{j}- c_{i})}{(1 -c_{i}a_{j})} = \prod_{\substack{1 \leq i \leq k\\l+1 \leq j \leq k}} (1 - c_{i}a_{j})\prod_{\substack{1 \leq j \leq l\\1 \leq i\leq k}}(a_{j}- c_{i}).$$
\end{enumerate}
It is probably useful to point out that in the above
(1) is $Z_{O}(C)$ and that (2) is $$\frac{\prod a_{j}^{k}}{Z_{S}(A -U +U^{-1},C)}$$ where $U = \{a_{1}, a_{2}, \cdots, a_{l} \}.$

All the rest of the factors only involve the terms with $a_{j}$'s and we turn to these now. 
First notice that 
\[ \frac{\prod_{\substack{1 \leq j \leq l\\1 \leq i\leq k}}(1 - a_{i}a_{j})}{\prod_{1\leq i\leq  j \leq k}( 1 - a_{i,}a_{j})\prod_{\substack{1 \leq i \leq l\\1 \leq j \leq l}}( 1/a_{j} - a_{i})}\]
yields when simplified 
\[ \frac{(-1)^{l^{2}}}{ \prod_{1}^{l}a_{j}\prod_{l+1 \leq i \leq j \leq k}( 1 - a_{i}a_{j}) \prod_{1 \leq i \leq j \leq l}( 1 - 1/a_{i}a_{j})} .\]
Next the remaining factors are
\[ \frac{\prod_{i = 2}^{l}\prod_{j = 1}^{i-1}( 1/a_{i} - 1/a_{j})( a_{i} - a_{j})}{\prod_{\substack{i\neq j\\ 1\leq j \leq l\\ 1 \leq i \leq k}}(a_{j} - a_{i})} .\]
This is the same as 
\[ \frac{\prod_{\substack{i \neq j\\ 1 \leq i \leq l\\ 1 \leq j \leq l}}( a_{i} - a_{j})}{\prod_{j} a_{j}^{l-1}\prod_{\substack{i\neq j\\ 1\leq j \leq l\\ 1 \leq i \leq k}}(a_{j} - a_{i})}\] or
 \[ \frac{1}{\prod_{j} a_{j}^{l-1}\prod_{\substack{1\leq j \leq l\\ l+1 \leq i \leq k}}(a_{j} - a_{i})}  = \frac{1}{\prod_{j} a_{j}^{l-1}a_{j}^{k - l}\prod_{\substack{ 1\leq j \leq l\\ l+1 \leq i \leq k}}(1 - a_{i}/a_{j})}.\]
 
 Putting this all together we have 
 \[ \frac{(-1)^{l^{2}}}{ \prod a_{j}^{k} \prod_{l+1 \leq i \leq j \leq k}( 1 - a_{i}a_{j}) \prod_{1 \leq i \leq j \leq l}( 1 - 1/a_{i}a_{j})\prod_{\substack{ 1\leq j \leq l\\ l+1 \leq i \leq k}}(1 - a_{i}/a_{j})},\]
 which is exactly 
 \[(-1)^{l}Z_{S}(A -U +U^{-1})/\prod a_{j}^{k} \]
 and this shows that for this set $U$ the contribution of the Fredholm determinant yields the same answer as in the ratios formula. 
 
 We have shown agreement with the coefficient of the term $a_{1}^{N}a_{2}^{N}\cdots a_{l}^{N}.$ The indexing is not important and hence this works for any subset $U$. Notice we have equated the $l \times l$ determinant in the Fredholm expansion with the sum over all subsets of size $l$ on the right-hand side of our formula in Theorem \ref{t1}.
 
 \endproof
 
 One final comment is that the in the proof we assumed that the number of factors in the numerator and denominator of $g$ were the same. But a moment's thought shows that if there are fewer or more $c_{i}$ terms than $a_{i}$ terms, the proof is easily modified to hold in those cases. 

This same proof with minor changes works in the orthogonal cases as well. Here is the statement in that case.

\begin{theorem}\label{t3} Let $A = \{a_{1}, a_{2}, \dots, a_{k}\}$ and let $C = \{c_{1}, a_{2}, \dots, c_{k}\}$ where $|c_i|  <1$. Let
 \[g(\theta) =  \prod_{i=1}^{k}\frac{(1 - a_{i}e^{i \theta})}{(1 - c_{i}e^{i \theta})}.\]  For a matrix $U$ with eigenvalues 
$\{e^{i\theta_j}:1\le j\le N\}$ we define 
$$\det g(U) = \prod_{1\le i \le k\atop 1\le j\le N} \frac{(1 - a_{i}e^{i \theta_j})}{(1 - c_{i}e^{i \theta_j})}
.$$Further, define $$Z_O(A;C)=\frac{Z_O(A)Z_S(C)}{Z_U(A,C)}.$$
\begin{enumerate}
\item Let $N = 2m +1.$ Then $$\int_{O^{+}(N)} \det g(U)
~dU=\sum_{U\subset A}(-1)^{|U|}U^N Z_O(A-U+U^{-1};C).$$
 \item Let $N = 2m +1.$ Then $$\int_{O^{-}(N)} \det g(U)
~dU=\sum_{U\subset A} U^N Z_O(A-U+U^{-1};C).$$
\item Let $N = 2m.$ Then $$\int_{O^{-}(N)} \det g(U)
~dU=\sum_{U\subset A}(-1)^{|U|}U^N Z_O(A-U+U^{-1};C).$$
\item Let $N = 2m.$ 
Then $$\int_{O^{+}(N)} \det g(U)
~dU=\sum_{U\subset A} U^N Z_O(A-U+U^{-1};C).$$
\end{enumerate}
\end{theorem}
We can do the same for the unitary case, although this is a bit harder. Before we begin the proof of the identity, we need three lemmas and some notation.

Let $\{a_{i}, a_{2}, \dots a_{l}\}$ and $ \{b_{i}, b_{2}, \dots b_{l}\}$ be any two ordered sets of complex numbers and define the $l \times l$ determinant by 
\[ D(\{a_{i}\}, \{b_{j}\}) = \det (a_{h}^{i_{g}} \,b_{h}^{i_{h}})_{ 1 \leq g, h \leq l}. \] 

Here is the $3 \times 3$ example for sets $\{a_{1}, a_{2}, a_{3}\}$ and $\{b_{1}, b_{2}, b_{3}\}$.

\[ D(\{a_{i}\}, \{b_{j}\}) = \left|\begin{array}{ccc} a_{1}^{i_{1}}b_{1}^{i_{1}} & a_{2}^{i_{1}}b_{2}^{i_{2}}& a_{3}^{i_{1}}b_{3}^{i_{3}}
 \\ a_{1}^{i_{2}}b_{1}^{i_{1}} & a_{2}^{i_{2}}b_{2}^{i_{2}}& a_{3}^{i_{2}}b_{3}^{i_{3}} \\
 a_{1}^{i_{3}}b_{1}^{i_{1}} & a_{2}^{i_{3}}b_{2}^{i_{2}}& a_{3}^{i_{3}}b_{3}^{i_{3}}\end{array}\right| .\]
 
 The reason for defining the above is that it will occur in our Fredholm expansion. 

\begin{lemma}
If any $a_{i} = a_{j}, i \neq j,$ then $D(\{a_{i}\}, \{b_{j}\})  = 0.$
\end{lemma}
This is easy to see since we can factor out $b_{h}^{i_{h}}$ from any column and the resulting determinant has two equal columns. 

\begin{lemma}
If any $b_{i} = b_{j}, i \neq j$ then $ \sum_{i_{i}, \dots ,i_{l} \geq m} D(\{a_{i}\}, \{b_{j}\})  = 0.$
\end{lemma}
 \proof
It suffices to assume that $b_{1} = b_{2}$. Factor first the $b_{h}^{i_{h}}$s from all the columns and then use the permutation definition of a determinant and consider the product of $a_{g}^{i_{h}}$ in each summand. There will be exactly one $a_{g}^{i_{h}}$ term with $i_{1}$ as an exponent and exactly one other $a_{g}^{i_{h}}$ term with $i_{2}$ as an exponent. This is because there can only be one element from each column. There is also a permutation with these exponents reversed and whose permutation sign is the opposite of the other. 
For this second permutation, if we change variables in the sum switching $i_{1}$ and $i_{2}$, we see that after summing the two terms cancel. 
\endproof
Notice the $D(\{a_{i}\}, \{b_{j}\})$ depends on the order of the $\{a_{i}\}, \{b_{j}\}$ in the determinant.
\begin{lemma}
Suppose there is a permutation of the  $\{b_{j}\}$s in the columns by the permutation $\sigma.$ Then 
$$ \sum_{i_{i}, \dots ,i_{l} \geq m} D(\{a_{i}\}, \{b_{j}\})$$ is changed by the sign of the permutation. 
\end{lemma}
\proof
Change both the $a_{j}$'s and $b_{j}$'s by the given permutation. After summing the answer does not change. Permuting the $a_{j} $s back just permutes the rows which changes the answer by the permutation sign and leaves the $b_{j}$'s still permuted. 
\endproof
\begin{theorem}\label{t2}
Let \[\phi = \prod_{i=1}^{k}\frac{(1 - a_{i}e^{-i \theta})(1 - b_{i}e^{i \theta})}{(1 - c_{i}e^{-i \theta})(1 - d_{i}e^{i \theta})},\] and let 
$A = \{a_{1}, a_{2}, \dots, a_{k}\}, B =  \{b_{1}, b_{2}, \dots, b_{k}\} ,C = \{c_{1}, c_{2}, \dots, c_{k}\}, D =  \{d_{1}, d_{2}, \dots, d_{k}\}$. Then 
$$ I_{U(N)}\phi =  \sum_{S\subset A ,T\subset B\atop |S|=|T|}S^N T^N    Z(A-S+T^{-1},B-T+S^{-1};C,D).
$$
\end{theorem}
\proof
We use the BOCG identity. First we note that by direct computation $E(\phi) = Z(A, B; C, D),$ which corresponds exactly to the empty set in the above. 

Referring to the identity the operator which occurs in the Fredholm determinant part of the identity is a product given by
$$ Q_{N}H(\phi_{-}\phi_{+}^{-1} )H(\tilde{\phi}_{-}^{-1} \tilde{\phi}_{+})Q_{N}.$$ The Fourier coefficients of $\tilde{\phi}_{-}^{-1} \tilde{\phi}_{+}$ are given by
\[ \sum_{j = 1}^{k}\al_{j}a_{j}^{l-1}, \,\,\,\,\,\, \al_{j} = \prod_{i \neq j}\frac{1}{(a_{j} - a_{i})}\prod_{i}\frac{(1 - b_{i}a_{j})(a_{j}- c_{i})}{(1 -d_{i}a_{j})}.\]
The Fourier coefficients of $\phi_{-}\phi_{+}^{-1}$ are given by
\[ \sum_{j = 1}^{k}\beta_{j}b_{j}^{l-1}, \,\,\,\,\,\, \beta_{j} = \prod_{i \neq j}\frac{1}{(b_{j} - b_{i})}\prod_{i}\frac{(1 - a_{i}b_{j})(b_{j}- d_{i})}{(1 -c_{i}b_{j})}.\]

Each $j,k,$ entry in $H(\tilde{\phi}_{-}^{-1} \tilde{\phi}_{+})$  is a sum of terms of the form $\al_{j}a_{j}^{j+k}$ and each $j,k$ entry in $H(\phi_{-}\phi_{+}^{-1} )$ is of the form  $\beta_{j}b_{j}^{j+k}.$ Thus when squaring the operator we find that the $j,k$ entry of the product is given by a sum of terms
\[ \frac{\beta_{g}\al_{h}b_{g}^{j}a_{h}^{k}}{(1-b_{g}a_{h})}.\]
Next consider the $l$th term in the Fredholm expansion,
$$\sum_{i_{1}, i_{2}, \cdots, i_{l}} \det K(i_{j}, i_{k})_{1 \leq j,k \leq l}.$$  Since $K$ is a sum of columns of the above form, we can expand $K(i_{j}, i_{k})$ in columns over all possible choices. Let us begin with the choice 
$\{a_{1}, a_{2}, \dots, a_{l}\}$ and $\{b_{1},b_{2}, \dots, b_{l}\},$ where the $b_{i}$ represents a subset of $B$.
By our lemmas we have no contribution if any of the $a_{i}$'s overlap or if any of the $b_{j}$'s overlap.
A moment's thought shows that $\det K(i_{j}, i_{k})$ is a sum of determinants, each of the form 
 $$  \frac{\prod_{j = 1}^{l}\al_{j}\beta_{j} }{\prod_{j=1}^{l}(1 - a_{j}b_{j})}D(\{b_{i}\}, \{a_{j}\}).$$
 We fix the order of the $b_{j}$'s. Then the contribution from all possible permutations of the $a_{j}$'s is 
 \[ \prod_{j = 1}^{l}\al_{j}\beta_{j} \,\,(\tilde{D}(\{b_{i}\} ,\{a_{j}\})\left(\sum_{\sigma}sgn(\sigma)\frac{1}{\prod_{j=1}^{l}(1 - b_{j}a_{\sigma(j)})}\right).\]
 Here $\tilde{D}$ represents the determinant $D(\{b_{i}\} ,\{a_{j}\})$ with the elements in the standard order. This expression on the right is a Cauchy determinant and after evaluating yields
  \[ \prod_{j = 1}^{l}\al_{j}\beta_{j} \,\,(\tilde{D}(\{b_{i}\},\{a_{j}\})\frac{\prod_{1 \leq i < j \leq l}(1/a_{i} - 1/a_{j})\prod_{1 \leq i < j \leq l}(b_{i} - b_{j})}{\prod_{i=1}^{i}a_{i}\prod_{i,j}(1/a_{i} - b_{j})}.\] After summing as before $\tilde{D}$ is also Cauchy and the result is that the above is
 \[ \prod_{j= 1}^{l}a_{j}^{N}b_{j}^{N}\prod_{j = 1}^{l}\al_{j}\beta_{j} \,\,\left(\frac{\prod_{1 \leq i < j \leq l}(1/a_{i} - 1/a_{j})\prod_{1 \leq i < j \leq l}(b_{i} - b_{j})}{\prod_{i=1}^{i}a_{i}\prod_{i,j}(1/a_{i} - b_{j})}\right)^{2}.\] Notice that this expression is independent of a permutation of the $b_{j}$'s. Hence we have $l!$ of these which cancels with $1/l!$ in the Fredholm expression.
 
 Using the notation of the theorem we think of the set $S$ as $\{a_{1}, a_{2}, \dots, a_{l}\}$ and $T$ as
 $\{b_{1}, b_{2}, \dots, b_{l}\}.$
 Now we just need to gather terms as before, and as before we have four types of terms. 
 \begin{enumerate}
 \item The terms with only $c_{j}$'s and $d_{j}$'s occurs only in $E(\phi)$ and is easily seen to be $Z(C,D).$
 \item The terms with both $a_{j}$'s, $b_{j}$'s and $d_{j}$'s is 
 $$ \prod_{i,j}(1 - a_{i}d_{j})\prod_{\substack{1 \leq j \leq l\\1 \leq i\leq k}}\frac{(b_{j}- d_{i})}{(1 -d_{i}a_{j})},$$
 which is easily seen to be 
 \[ \frac{\prod_{i = 1}^{l}b_{i}^{k} }{Z( A - S + T^{-1}, D)}.\]
 \item The terms with both $a_{j}$'s, $b_{j}$'s and $c_{j}$'s is 
 $$ \prod_{i,j}(1 - b_{i}c_{j})\prod_{\substack{1 \leq j \leq l\\1 \leq i\leq k}}\frac{(a_{j}- c_{i})}{(1 -b_{i}c_{j})},$$
 or \[ \frac{\prod_{i = 1}^{l}a_{i}^{k} }{Z( B - T + S^{-1}, D)}.\]
   \end{enumerate}
 Hence the only thing left to consider is the terms with only $a_{j}$'s and $b_{j}$'s. This is 
 \[ Z(A,B) \prod_{\substack{i\neq j\\ 1\leq j \leq l\\ 1 \leq i \leq k}}\frac{1}{(a_{j} - a_{i})}\prod_{\substack{1 \leq j \leq l\\1 \leq i\leq k}}(1 - b_{i}a_{j}) \prod_{\substack{i\neq j\\ 1\leq j \leq l\\ 1 \leq i \leq k}}\frac{1}{(b_{j} - b_{i})}\prod_{\substack{1 \leq j \leq l\\1 \leq i\leq k}}(1 - b_{j}a_{i})\\\]
\[ \times \left(\frac{\prod_{1 \leq i < j \leq l}(1/a_{i} - 1/a_{j})\prod_{1 \leq i < j \leq l}(b_{i} - b_{j})}{\prod_{i=1}^{i}a_{i}\prod_{i,j}(1/a_{i} - b_{j})}\right)^{2}.\]
After simplifying the above, we find that it becomes
\[ (-1)^{l}Z(A - S + T^{-1}, B -T +S^{-1})/\prod_{i = 1}^{l} (a_{i}b_{i})^{k} \] and after combining with the other terms yields the correct answer. 
\endproof

\section{Appendix}
Here we give some additional computations.
  
As the simplest interesting example of a ratios factorization, suppose that $A=\{a\},B=\{b\},C=\{c\},D=\{d\}$  where $|c|,|d|<1.$ By Theorem 6.6,
\begin{eqnarray*}
\int_{U(N)}\frac{  \Lambda_U(a)  \Lambda_{U^*}(b)}
{ \Lambda_U(c)  \Lambda_{U^*}(d)}dU 
&=&\frac{(1-ad)(1-bc)}{(1-ab)(1-cd)}+a^Nb^N\frac{(1-\frac db)(1-\frac ca)}{(1-\frac1{ab})(1-cd)}.
\end{eqnarray*}
 We are interested in the case that $a=b$ and $c=d$; we assume this from now on. 
 We have
 \begin{eqnarray*}
 R_{U(N)}(a;c)&:=&\int_{U(N)}\frac{  \Lambda_U(a)  \Lambda_{U^*}(a)}
{ \Lambda_U(c)  \Lambda_{U^*}(c)}dU\\
& =& \frac{(1 + a^{(1 + N)} - a c - a^N c) (1 - a^{(1 + N)} - a c + a^N c)}{(1 - a) (1 + a) (1 - c) (1 + c)}
.\end{eqnarray*}

Thus,
\begin{eqnarray*}
\begin{array}{|c|l|} 
\hline
N & R_U(N) \\
\hline
1& \frac{1-2ac+a^2}{1-c^2}\\
\hline
2&\frac{(1-a+a^2-ac)(1+a+a^2-ac)}{1-c^2}\\
\hline
3&\frac{(1 + a^2 - a c) (1 + a^4 - a c - a^3 c)}{1-c^2}\\
\hline
4&\frac{(1 + a + a^2 + a^3 + a^4 - a c - a^2 c - a^3 c) (1 - a + a^2 - a^3 + 
   a^4 - a c + a^2 c - a^3 c)}{1-c^2}\\
   \hline
\end{array}
\end{eqnarray*}

 For the symplectic group we have by Theorem 6.1,
 \begin{eqnarray*}
R_{S(N)} &:=&  \int_{USp(N)}\frac{\Lambda_U(a)}{\Lambda_U(c)}~dU\\
 &=& \frac{(1-ac)}{(1-a^2) }+a^N \frac{(1-c/a)}{(1-a^{-2}) }.
 \end{eqnarray*}

Thus
\begin{eqnarray*}
\begin{array}{|c|l|}\hline
N&R_S(N)\\
\hline
2&1  + a^2 - a c \\
\hline
4& 1+a^2+a^4-ac-a^3 c\\
\hline
6&  1+a^2+a^4-+a^6-ac-a^3 c-a^5 c\\
\hline
8& 1+a^2+a^4+a^6+a^8-ac-a^3 c-a^5 c-a^7 c\\
\hline
\end{array}
\end{eqnarray*}

Using Theorem \ref{t3} and an obvious notation for the orthogonal averages we have
\begin{eqnarray*}
\begin{array}{|c|l|}
\hline
N& R_{O^{+}(N)}\\
\hline
1& \frac{1 - a }{1-c}\\
\hline
2& \frac{1 + a^2 - 2a c }{1-c^2}\\
\hline
3& \frac{1 -ac -a^{3} +a^{2}c}{1-c^2}\\
\hline
4& \frac{1   - a c +a^{4}- a^3 c}{1-c^2}\\
\hline
5&\frac{1 -ac -a^{5} +a^{4}c}{1-c^2}\\
\hline
6& \frac{1 -ac +a^{6} -a^{5}c}{1-c^2}\\
\hline
\end{array}
\end{eqnarray*}
and 
\begin{eqnarray*}
\begin{array}{|c|l|}
\hline
N& R_{O^{-}(N)}\\
\hline
1& \frac{1 + a }{1+c}\\
\hline
2& \frac{1 -a^2 }{1-c^2}\\
\hline
3& \frac{1 -ac +a^{3} -a^{2}c}{1-c^2}\\
\hline
4& \frac{1   - a c -a^{4}+a^3 c}{1-c^2}\\
\hline
5&\frac{1 -ac +a^{5} -a^{4}c}{1-c^2}\\
\hline
6& \frac{1 -ac -a^{6} +a^{5}c}{1-c^2}\\
\hline
\end{array}
\end{eqnarray*}

 From these we can see the verification of our factorization formulas.

\begin{eqnarray*}R_{U(1)}&=& R_{S(0)}R_{O^{+}(2)}\\
R_{U(2)}&=& (g(1)g(-1))^{-1}R_{O^{+}(3)}R_{O^{-}(3)}\\
R_{U(3)}&=& R_{S(2)}R_{O^{+}(4)}\\
R_{U(4)}&=& (g(1)g(-1))^{-1}R_{O^{+}(5)}R_{O^{-}(5)}\\
R_{U(5)}&=& R_{S(4)}R_{O^{+}(6)}.\\
\end{eqnarray*}

Here  $g(\pm 1)$  means substitute $\pm 1$  for $e^{i\theta}$ in the formula for $g$. So for
example, if
$g = \frac{1 - a e^{i\theta}}{1-ce^{i\theta}}$  then $g(1) = \frac{1 - a}{1-c}$.

We now consider the slightly more complicated case of  unitary moments with two numerator parameters $a$ and $b$
$$I_{U(N)}(a,b):= \int_{U(M)} \Lambda_U(a) \Lambda_U(b) \Lambda_{U^*}(a) \Lambda_{U^*}(b)~dU.$$
By Theorem 6.6 with $A=B=\{a,b\}$ and $C=D=\{\}$, the empty set, we have
  \begin{eqnarray*}
   I_{U(1)}(a,b)&=&1 + a^2 + 2 a b + b^2 + a^2 b^2\\
I_{U(2)}(a,b)&  =&(1 - a + a^2 - b + 2 a b - a^2 b + b^2 - a b^2 + a^2 b^2)\\
& &\quad \times (1 + a + 
   a^2 + b + 2 a b + a^2 b + b^2 + a b^2 + a^2 b^2)\\
I_{U(3)}(a,b)&=&(1 + a^2 + a b + b^2 + a^2 b^2)\\
& &\quad \times   (1 + a^4 + a b + a^3 b + 2 a^2 b^2 + 
   a b^3 + a^3 b^3 + b^4 + a^4 b^4)\\
   I_{U(4)}(a,b)    &= &(1 - a + a^2 - a^3 + a^4 - b + 2 a b - 2 a^2 b
 + 2 a^3 b - a^4 b\\ 
 & &\qquad + 
   b^2 - 2 a b^2 + 3 a^2 b^2 - 2 a^3 b^2 + a^4 b^2 
    b^3+ 2 a b^3 - 
   2 a^2 b^3 \\
   & &\qquad + 2 a^3 b^3 - a^4 b^3 + b^4 - a b^4 + a^2 b^4 - a^3 b^4 +
    a^4 b^4)\\
    & &\quad \times  (1 + a + a^2 + a^3 + a^4 + b + 2 a b + 2 a^2 b + 
   2 a^3 b \\
   & &\qquad + a^4 b + b^2 + 2 a b^2 + 3 a^2 b^2  
+ 2 a^3 b^2 + a^4 b^2 +
    b^3\\
    & & \qquad + 2 a b^3 + 2 a^2 b^3 + 2 a^3 b^3 + a^4 b^3 + b^4 + a b^4 + 
   a^2 b^4 \\
   &&\qquad + a^3 b^4 + a^4 b^4).
   \end{eqnarray*}

   For the symplectic moments we have   

 \begin{eqnarray*}I_{USp(2)}(a,b)&=&1 + a^2 + a b + b^2 + a^2 b^2\\
I_{USp(4)}(a,b)&=&1 + a^2 + a^4 + a b + a^3 b + b^2 + 2 a^2 b^2 + a^4 b^2 + a b^3\\&&\qquad  + 
 a^3 b^3 + b^4 + a^2 b^4 + a^4 b^4.\\
 \end{eqnarray*}
 
 For the orthogonal moment we have, 
 \begin{eqnarray*}
I_{O^{-}(3)}(a,b)& = & (1 + a) (1 + b) (1 - a + a^2 - b + 2 a b - a^2 b\\&&\qquad + b^2 - a b^2 + 
   a^2 b^2)\\
I_{O^{+}(3)}(a,b)&=& (1 - a) (1 - b) (1 + a + a^2  +b + 2 a b + a^2 b\\&&\qquad  + b^2 + a b^2 + 
   a^2 b^2)\\
I_{O^{-}(5)}(a,b)&=&
(1+a)(1+b)(1 - a + a^2 - a^3 + a^4 - b + 2 a b \\&&\qquad - 2 a^2 b
 + 2 a^3 b - a^4 b + 
   b^2  - 2 a b^2 + 3 a^2 b^2 \\&&\qquad - 2 a^3 b^2 + a^4 b^2 
    b^3+ 2 a b^3 - 
   2 a^2 b^3 +
   2 a^3 b^3  \\&&\qquad \quad- a^4 b^3 + b^4 - a b^4 + a^2 b^4    - a^3 b^4 +
    a^4 b^4) \\
    I_{O^{+}(5)}(a,b)&=&(1 - a) (1 - b) (1 + a + a^2 + a^3 + a^4 + b + 2 a b\\&&\qquad  + 2 a^2 b + 
   2 a^3 b 
   + a^4 b + b^2 + 2 a b^2 + 3 a^2 b^2  \\&&\qquad
+ 2 a^3 b^2 + a^4 b^2 +
    b^3
    + 2 a b^3+ 2 a^2 b^3  \\&&\qquad + 2 a^3 b^3 + a^4 b^3  + b^4 + a b^4 + 
   a^2 b^4 + a^3 b^4 + a^4 b^4) \\
    I_{O^{+}(4)}(a,b)&=&1 + a^4 + a b + a^3 b + 2 a^2 b^2 + a b^3 + a^3 b^3 + b^4 + a^4 b^4.\end{eqnarray*}

Thus, we have the following identities:
\begin{eqnarray*}I_{U(1)}(a,b)&=&I_{O^{+}(2)}(a,b)\\
  I_{U(2)}(a,b)&=& (g(1)g(-1))^{-1}I_{O^{-}(3)}(a,b)I_{O^{+}(3)}(a,b)\\
  I_{U(3)}(a,b)&=&I_{USp(2)}(a,b)I_{O^{+}(4)}(a,b)\\
  I_{U(4)}(a,b)&=&(g(1)g(-1))^{-1}I_{O^{-}(5))}(a,b)I_{O^{+}(5)}(a,b).
  \end{eqnarray*}

 \end{document}